\documentclass{aims}
\usepackage{amsmath}
 \usepackage{paralist}
 \usepackage{graphics} %% add this and next lines if pictures should be in esp format
 \usepackage{epsfig} %For pictures: screened artwork should be set up with an 85 or 100 line screen
\usepackage[colorlinks=true]{hyperref}
\hypersetup{urlcolor=blue, citecolor=red}

\usepackage{epsfig}
\usepackage{graphicx}
\usepackage{amsmath}
\usepackage{amssymb}
\usepackage{setspace}
\usepackage{graphicx}
\graphicspath{{/}}
\usepackage{graphicx}

\newtheorem{theorem}{Theorem}[section]

\theoremstyle{definition}

\title[Conformal Approach for Surface Inpainting]
     {A Conformal Approach for Surface Inpainting}

\author[Lui, Wen, Gu]{\small{$\textsc{Lok\ Ming\ Lui}^{1}$, \small{$\textsc{ Chengfeng\ Wen}^1$} and \small{$\textsc{Xianfeng\ Gu}^{2}$}}\\
\textit{\footnotesize{$^1$ Department of Mathematics, The Chinese University of Hong Kong}}\\
\textit{\footnotesize{$^2$ Department of Computer Science, State University of New York at Stony Brook}}\\
}

\begin{document}
\maketitle

% Enter the first author's name and address:
%\centerline{\scshape First-name1 last-name1 }
%\medskip
%{\footnotesize
% please put the address of the first author
% \centerline{First line of the address of the first author}
%   \centerline{Other lines}
%   \centerline{ Springfield, MO 65801-2604, USA}
%} % Do not forget to end the {\footnotesize by the sign }

%\medskip

%\centerline{\scshape First-name2 last-name2 and First-name3
%last-name3}
%\medskip
%{\footnotesize
% please put the address of the second  and third author
% \centerline{ First line of the address of the second author}
%  \centerline{Other lines}
%  \centerline{Springfield, MO 65810, USA}
%}

%\bigskip

% The name of the associate editor will be entered by an editorial staff
% \centerline{(Communicated by the associate editor name)}

%The abstract of your paper
\begin{abstract}
We address the problem of surface inpainting, which aims to fill in holes or missing regions on a Riemann surface based on its surface geometry. In practical situation, surfaces obtained from range scanners often have holes or missing regions where the 3D models are incomplete. In order to analyze the 3D shapes effectively, restoring the incomplete shape by filling in the surface holes is necessary. In this paper, we propose a novel conformal approach to inpaint surface holes on a Riemann surface based on its surface geometry. The basic idea is to represent the Riemann surface using its conformal factor and mean curvature. According to Riemann surface theory, a Riemann surface can be uniquely determined by its conformal factor and mean curvature up to a rigid motion. Given a Riemann surface $S$, its mean curvature $H$ and conformal factor $\lambda$ can be computed easily through its conformal parameterization. Conversely, given $\lambda$ and $H$, a Riemann surface can be uniquely reconstructed by solving the Gauss-Codazzi equation on the conformal parameter domain. Hence, the conformal factor and the mean curvature are two geometric quantities fully describing the surface. With this $\lambda$-$H$ representation of the surface, the problem of surface inpainting can be reduced to the problem of image inpainting of $\lambda$ and $H$ on the conformal parameter domain. The inpainting of $\lambda$ and $H$ can be done by conventional image inpainting models. Once $\lambda$ and $H$ are inpainted, a Riemann surface can be reconstructed which effectively restores the 3D surface with missing holes. Since the inpainting model is based on the geometric quantities $\lambda$ and $H$, the restored surface follows the surface geometric pattern as much as possible. We test the proposed algorithm on synthetic data, 3D human face data and MRI-derived brain surfaces. Experimental results show that our proposed method is an effective surface inpainting algorithm to fill in surface holes on an incomplete 3D models based their surface geometry.
\end{abstract}

%Introduction
\section{Introduction} Surface processing is an important topic in computer graphics, computer visions and medical imaging. Examples includes surface denoising, surface inpainting, surface remeshing and so on. Amongst the various tasks in surface processing, the problem of surface inpainting, which aims to fill in holes or missing regions on a Riemann surface, is an important pre-processing step necessary for the geometric analysis of incomplete 3D shapes. Digital 3D surfaces containing holes are common. In practical situation, real-life objects are usually captured from range scanners. Despite the recent advance in the laser scanning technology, 3D models obtained from even the most sophisticated acquisition devices are often incomplete due to occlusions, low reflectance, constraints during recording or self-occlusions of complicated geometry. This results in holes or missing regions in the obtained 3D digital models. In order to analyze the geometric objects effectively or derive a visually appealing 3D models, developing an effective algorithm for filling in holes or missing regions on a surface is of utmost importance.

In this paper, we propose a framework for filling in holes or missing regions on a Riemann surface using conformal geometry. The ultimate goal is to obtain a restored surface with holes filled based on the surface geometry, so that the basic geometry of the surface can be smoothly patched. For this purpose, a good representation of the Riemann surface describing its geometry is required. In this work, we propose to represent the Riemann surface using its conformal factor and mean curvature. According to Riemann surface theory, a Riemann surface can be uniquely determined by its conformal factor and mean curvature up to a rigid motion. Given a Riemann surface $S$, its mean curvature $H$ and conformal factor $\lambda$ can be computed easily through its conformal parameterization. Conversely, given $\lambda$ and $H$, a Riemann surface can be uniquely reconstructed by solving the Gauss-Codazzi equation on the conformal parameter domain. Hence, the conformal factor and the mean curvature are two geometric quantities fully describing the surface. Given the $\lambda$-$H$ representation of an incomplete surface, the problem of surface inpainting can be reformulated into the image inpainting problem of $\lambda$ and $H$. Using the conventional image inpainting techniques, both $\lambda$ and $H$ can be inpainted on the conformal parameter domain. A Riemann surface can then be reconstructed which restores the 3D surface with missing holes. Since the proposed inpainting model is based on two geometric quantities fully determining the surface, the surface holes can be effectively filled , which follows the surface geometric pattern. We test the proposed algorithm on synthetic data, 3D human face data and MRI-derived brain surfaces. Experimental results show that our algorithm can effectively inpaint surface holes based on the surface geometry to restore the incomplete 3D surface models.

The contributions of this paper are two-folded. First, we consider the representation of a Riemann surface using its conformal factor $\lambda$ and mean curvature $H$. Given two scalar functions $\lambda$ and $H$, the associated Riemann surface $S$ can be reconstructed through solving the Gauss-Codazzi equation numerically. Second, we propose a surface inpainting algorithm by inpainting the scalar functions $\lambda$ and $H$. The surface inpainting algorithm is based on $\lambda$ and $H$, which are two geometric quantities fully describing the surface. Hence, the proposed surface inpainting algorithm can fill in surface holes based on the surface geometry.

Our paper is organized as follows: prior work on related topics is presented in section 2. The basic mathematical theory is discussed in section 3. In Section 4, the details of our proposed surface representation using conformal factor and mean curvature are discussed. The detailed algorithm for surface inpainting is described in Section 5. In Section 6, the numerical algorithms are summarized. Experimental results are reported in Section 7, and some conclusions and future work are discussed in Section 8.

\section{Previous work}
In this section, we will review some relevant works closely related to our paper.

The problem of surface inpainting was inspired by 2D image completion, which has been widely studied by different groups. Different variational models have been proposed \cite{Chaninpainting,Chaninpainting2,ChanHO,ChanHO2,MasnouHO}. The basic idea is to diffuse the image intensity into the inpainting domain from its surrounding regions based on various regularizations of the solution. For example, Chan et al. \cite{Chaninpainting,Chaninpainting2} proposed the total variation (TV) models for image inpainting. Euler elastica model \cite{ChanHO,ChanHO2,MasnouHO} has also been proposed, which is essentially a variational model with a regularizing term containing higher-order derivatives. Other examples include the active contour model based on Mumford and Shahs segmentation \cite{MS}, the inpainting scheme based on the Mumford-Shah-Euler image model \cite{Esedoglu}, the inpainting model with the Navier-Stokes equation \cite{NS}, and wavelet-based inpainting algorithms \cite{Waveletinpainting,Waveletinpainting2}. For a detailed introduction to image inpainting, we refer the readers to \cite{ChanHO}.

3D surface processing has also been extensively studied. Different surface denoising algorithms have been proposed. For example, surface
smoothing by mean curvature flow or Laplace smoothing have been proposed \cite{MCF1,MCF2,MCF3}, which are effective for removing oscillations on noisy surfaces. In order to preserve ridges and sharp corners of the surfaces, different second order anisotropic surface diffusion models have been considered \cite{2order,2order2,2order3,2order4}, which minimizes the weighted surface area. Higher-order isotropic or anisotropic flows based on minimizing the surface total curvature or weighted total curvature have also been discussed \cite{4order1,4order2}. Besides, various surface inpainting algorithms to fill in surface holes have been proposed \cite{InpaintingCaselles,InpaintingClarenz,InpaintingDavis,InpaintingSharf,InpaintingSavchenko,InpaintingVerdera}. A common approach is to fill in the occluded regions with surface patches, which smoothly attach to the boundary vertices of the surface holes \cite{InpaintingClarenz,InpaintingDavis}. For example, Clarenz et al. \cite{InpaintingClarenz} proposed to fill in surface holes by minimizing the Willmore energy functional. The surface patches obtained is guaranteed to satisfy the continuity properties. Davis et al. \cite{InpaintingDavis} proposed to extend a signed distance function, which is initially defined only on domains close to the known surface, to the complete space using volumetric diffusion. The surface can then be completed, even for non-trivial hole boundaries. Later, Caselles et al. \cite{InpaintingCaselles} proposed a level set variational framework based on the minimization of the $L^p$ norm of the mean curvature for surface inpainting. The algorithm produces an inpainted surface with minimal surface area given the boundary constraints, although the restored surface may produce non-smooth global geometry. Most of the above algorithms do not incorporate or only partially incorporate the geometric information of the surface into their models. The inpainting result obtained usually does not follow surface geometric pattern, and thus the restored surface might often be unnatural. In this paper, we propose to inpaint surface holes through inpainting the conformal factor and mean curvature, which are two geometric quantities fully determining the surface. The inpainted surface can better follow the surface geometric pattern.

The proposed surface inpainting model in this paper is based on the conformal structure of the Riemann surface. The problem of computing the conformal structure of a Riemann surface has been extensively studied, and different algorithms have been proposed \cite{Fischl2, Gu1, Gu3, Haker, Hurdal, Gu2}. For example, Hurdal et al. \cite{Hurdal} proposed to compute the conformal parameterization using circle packing and applied it to the registration of human brains. Gu et al. \cite{Gu1, Gu3, Gu2} proposed to compute the conformal parameterizations of human brain surfaces for registration using harmonic energy minimization and holomorphic 1-forms. Using conformal factor and curvatures, a shape index has also been proposed to measure geometric difference between hippocampal surfaces \cite{LuiBHF,LuiBHFHP}. Conformal structure has further been applied for geometric compression. In \cite{Gucompression}, the author proposed to compute the conformal structure using holomorphic one-form. The conformal structure is then applied to estimate the edge lengths of the triangulation mesh. Together with the angle between triangular faces, the surface mesh can be reconstructed by solving a linear system. Hence, instead of 3 coordinate functions, we can store the surface using only the conformal structure and angles between triangular faces.

\section{Mathematical background} In this section, we describe some basic mathematical concepts relevant to describing our algorithm.

\subsection{Conformal structure of a Riemann surface}
A surface $S$ with a conformal structure is called a {\it Riemann surface}. All Riemann surfaces are locally Euclidean. Given two Riemann surfaces $M$ and $N$. We can represent them locally as:

\begin{equation}
\phi_M(x_1,x_2): R^2\to M\subseteq \mathbb{R}^3 \mathrm{\ and\ } \phi_N(x_1, x_2): R^2\to N\subseteq \mathbb{R}^3
\end{equation}

\noindent respectively, where $(x_1,x_2)$ are their coordinates. The inner product of the tangent vectors at each point of the surface can be represented by its first fundamental form. The first fundamental form on $M$ can be written as
\begin{equation}
ds_M^2=\sum_{i,j} g_{ij} dx^i dx^j,
\end{equation}

\noindent where $g_{ij}=\frac{\partial \phi_M}{\partial x^i}\cdot \frac{\partial \phi_M}{\partial x^j}$ and $i,j=1,2$. Similarly, the first fundamental form on $N$ can be written as
\begin{equation}
ds_N^2=\sum_{i,j} \widetilde{g}_{ij} dx^i dx^j,
\end{equation}
\noindent where $g_{ij}=\frac{\partial \phi_N}{\partial x^i}\cdot \frac{\partial \phi_N}{\partial x^j}$ and $i,j=1,2$.

Given a map $f:M\to N$ between the $M$ and $N$. With the local parameterization, $f$ can be represented locally by its coordinates as $f:R^2\to R^2$, $f(x_1,x_2)=(f_1(x_1 , x_2), f_2(x_1,x_2))$. Every tangent vectors $\vec{v}$ on $M$ can be mapped ($push\ forward$) by $f$ to a tangent
vectors $f_*(\vec{v})$ on $N$. The inner product of the vectors $f_*(\vec{v_1})$ and $f_*( \vec{v_2}))$, where $\vec{v_1}$ and
$\vec{v_2}$ are tangent vectors on $M$, is:

\begin{equation}
\begin{split}
f^*(ds_N^2)(v_1,v_2) & :=\ <f_*(v_1), f_*(v_2)> \\
& =  \sum_{i,j} \widetilde{g}_{ij} f_*(v_i)\cdot f_*(v_j) \\
& =  \sum_{i,j} (\sum_{m,n}\widetilde{g}_{mn}\frac{\partial
f_i}{\partial x^m}\frac{\partial f_j}{\partial x^n}) v_i v_j)
\end{split}
\end{equation}

Therefore, a new Riemannian metric $f^*(ds_N^2)$ on $M$ is induced by $f$ and $ds_N^2$, called the $pull\ back\ metric$. We say that
the map $f$ is {\it conformal} if

\begin{equation}
f^*(ds_N^2)=\lambda(x_1, x_2)^2ds_M^2
\end{equation}

$\lambda$ is called the {\it conformal factor}. It is one important geometric quantity describing the Riemann surface $S$.

A parameterization $\varphi: R^2\to M$ is a {\it conformal parameterization} if $\varphi$ is a conformal map.

\begin{figure}[t]
\centering
\includegraphics[height=1.65in]{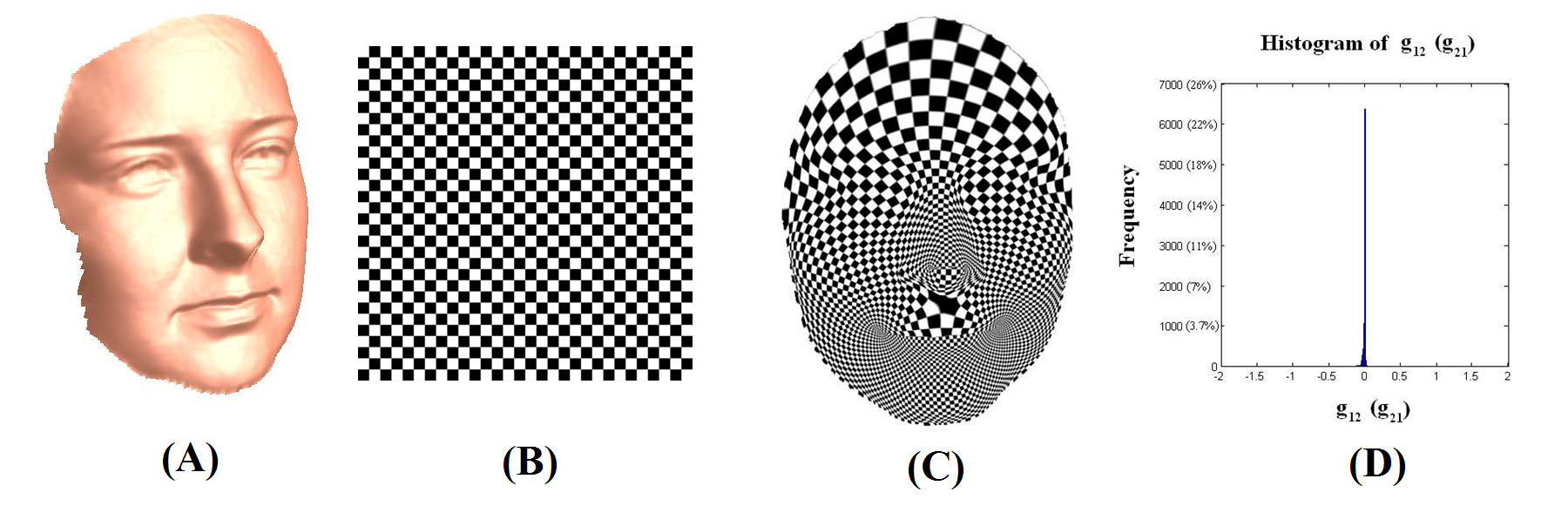}
\caption{Example of conformal parameterization. (A) shows a human face. It is conformally parameterized onto the 2D rectangle as shown in (B). We map the texture in (B) onto the surface, as shown in (C). Note that the right angle structure of the checkerboard texture is well-preserved. (D) shows the histogram of $g_{12}$ (or $g_{21}$) of the Riemannian metric. \label{fig:Confexample}}
\end{figure}

Intuitively, a map is conformal if it preserves the inner product of the tangent vectors up to a scaling factor, called the conformal
factor, $\lambda$. An immediate consequence is that every conformal map preserves angles. Figure \ref{fig:Confexample} gives an example of conformal parameterization of a human face.  (A) shows a human face. It is conformally parameterized onto the 2D rectangle as shown in (B). We map the texture in (B) onto the surface, as shown in (C). Note that the right angle structure of the checkerboard texture is well-preserved. (D) shows the histogram of $g_{12}$ (or $g_{21}$) of the Riemannian metric. It is concentrated at 0, meaning that the mapping is conformal.

\subsection{Curvatures of a Riemann surface}
We next briefly describe the concept of curvatures on a Riemann surface. Given a Riemann surface $S$, the normal curvature $\kappa_n$ in some direction is the reciprocal of the radius of the circle that best approximate a normal slice of the surface in that direction. For smooth surface, it can be computed from a $2\times 2$ symmetric matrix $\mathbf{W}$, called the {\it Weingarten matrix}, by:

\begin{equation}
\kappa_n = (u,v)\mathbf{W}\left(\begin{array}{c}
u\\
v
\end{array}
\right) = (u,v)\left(\begin{array}{cc}
e & f\\
f& g\\
\end{array}\right)\left(\begin{array}{c}
u\\
v
\end{array}
\right)
\end{equation}

\noindent for any unit length vector $(u, v)$ in the tangent plane of the surface.

The principal curvatures, principal directions, mean curvatures and Gaussian curvatures can be defined by the eigenvalues and eigenvectors of $\mathbf{W}$. Let $k_1, k_2$ be the eigenvalues of $\mathbf{W}$ and $\vec{d}_1$, $\vec{d}_2$ be the eigenvectors of $\mathbf{W}$. Then, $k_1$ and $k_2$ are called the {\it principal curvatures}; $\vec{d}_1$ and $\vec{d}_2$ are called the {\it principal directions}. The {\it mean
curvature} $H$ is defined as the average of the principal curvatures:
\begin{equation}
H = \frac{k_1+k_2}{2}
\end{equation}

The {\it Gaussian curvature} $K$ is defined as the product of the principal curvatures:
\begin{equation}
K = k_1 k_2
\end{equation}

With the conformal parameterization $\phi:D\subset \mathbb{R}^2\to S$ of $S$, the Gaussian curvature $K$ and the mean curvature $H$ can be computed easily:
\begin{equation}
K = -\frac{1}{2\lambda^2} \Delta \log \lambda
\end{equation}
\noindent where $\lambda$ is the conformal factor. And

\begin{equation}\label{H}
H = \frac{1}{2\lambda^2} \mathbf{sign}(\phi) |\Delta \phi|
\end{equation}

\noindent where $\mathbf{sign}(\phi) = \frac{<\Delta \phi, \vec{n}>}{|\Delta \phi|}$ and $\vec{n}$ is the (unit) surface normal.

It turns out the mean curvature $H$ together with the conformal factor $\lambda$ are two geometric quantities fully describing the surface $S$.

\section{Surface representation by conformal factor and mean curvature}
In this section, we will describe how a Riemann surface $S$ can be represented by its conformal factor $\lambda$ and mean curvature $H$.

Every Riemann surface is locally Euclidean, and hence it can be locally parametrized. Under the parameterization, a Riemann surface $S$ embedded in $\mathbb{R}^3$ can be represented by the position vector $\phi(u,v)$:

\begin{equation}
\phi(u,v) = (X(u,v), Y(u,v), Z(u,v)) \in \mathbb{R}^3.
\end{equation}

According to the Riemann surface theory, the surface $S$ can be determined uniquely (up to a rigid motion) by its first fundamental form $ds^2$ and second fundamental form.  The first fundamental form can be written as follows:

\begin{equation} \label{I}
ds^2 = E(u,v) du^2 + 2F(u,v) du dv + G(u,v) dv^2
\end{equation}

\noindent where $E(u,v) = <\phi_u , \phi_u>$, $F(u,v) = <\phi_u , \phi_v>$ and $E(u,v) = <\phi_v , \phi_v>$.

If $\phi$ is a conformal parameterization, $F(u,v) = 0$ and $E(u,v) = F(u,v) = \lambda(u,v)$, where $\lambda$ is the conformal factor. Hence,

\begin{equation}
ds^2 = \lambda(u,v)^2 (du^2 + dv^2)
\end{equation}

The second fundamental form can be written as follows:

\begin{equation}
II = L(u,v) du^2 + 2M(u,v) du dv + N(u,v) dv^2
\end{equation}

\noindent where $L(u,v) = <\phi_{uu} , \vec{n}>$, $M(u,v) = <\phi_{uv} , \vec{n}>$ and $N(u,v) = <\phi_{vv} , \vec{n}>$ ($\vec{n}$ is the unit surface normal).

Now, our goal is to find geometric quantities which can determine the position vector $\phi:=(X,Y,Z)\in \mathbb{R}^3$.

Let $z = u + i v$ ($i=\sqrt{-1}$), $dz = du + i dv$, $d\bar{z} = du - i dv$, $\frac{\partial}{\partial z}  = \frac{1}{2}(\frac{\partial}{\partial u} - i \frac{\partial}{\partial v})$ and $\frac{\partial}{\partial \bar{z}}  = \frac{1}{2}(\frac{\partial}{\partial u} + i \frac{\partial}{\partial v})$.

Suppose $\phi$ is a conformal parameterization, then we must have the following:

\begin{equation}
\begin{split}
<\phi_z, \phi_z> &= 0\\
<\phi_{\bar{z}}, \phi_{\bar{z}}> &= 0\\
<\phi_z, \phi_{\bar{z}}> &= \frac{\lambda^2}{2}.
\end{split}
\end{equation}

Besides, the natural frame $(\phi_z, \phi_{\bar{z}}, \vec{n})$ must satisfy certain equations. Consider $\mu = <\phi_{zz},\vec{n}>$. Then, the natural frame satisfies:
\begin{equation}\label{frame}
\frac{\partial}{\partial z}\left(\begin{array}{c}
\phi_z\\
\phi_{\bar{z}}\\
\vec{n}
\end{array}
\right) = \left(\begin{array}{c}
\frac{2}{\lambda}\lambda_z \phi_z  + \mu \vec{n}\\
\frac{\lambda^2}{2} H \vec{n} \\
- H \phi_z - \frac{2\mu \phi_{\bar{z}}}{\lambda^2}
\end{array}
\right)
\end{equation}

Clearly, the natural frame can be determined by the conformal factor $\lambda$, mean curvature $H$ and $\mu$. Besides, $\mu$ can as well be determined by $\lambda$ and $H$ through the Codazzi equation:

\begin{equation}\label{Codazzi}
\mu_{\bar{z}} = \frac{\lambda^2}{2} H_z
\end{equation}

Differentiate equation (\ref{Codazzi}) with respect to $z$, we obtain

\begin{equation}\label{Codazzi2}
\Delta \mu = \frac{1}{2}\lambda (2 \lambda_z H_z + \lambda H_{zz}).
\end{equation}

Given the value of $\phi$ on the boundary neighborhood $U$, the natural frame $(\phi_z, \phi_{\bar{z}}, \vec{n})$ can be found by solving the system of partial differential equation (\ref{frame}). The position vector $\phi$ can be reconstructed by integrating the natural frame. Therefore, the mean curvature $H$ and conformal factor $\lambda$ are two geometric quantities fully determining a Riemann surface. Given a Riemann surface $S$, its associated ($\lambda, H$) can be computed. Conversely, given ($\lambda, H$), the associated Riemann surface $S$ can be reconstructed. We call the pair ($\lambda, H$) the {\it  $\lambda$-$H$ representation} of $S$.

\begin{figure}[t]
\centering
\includegraphics[height=1.85in]{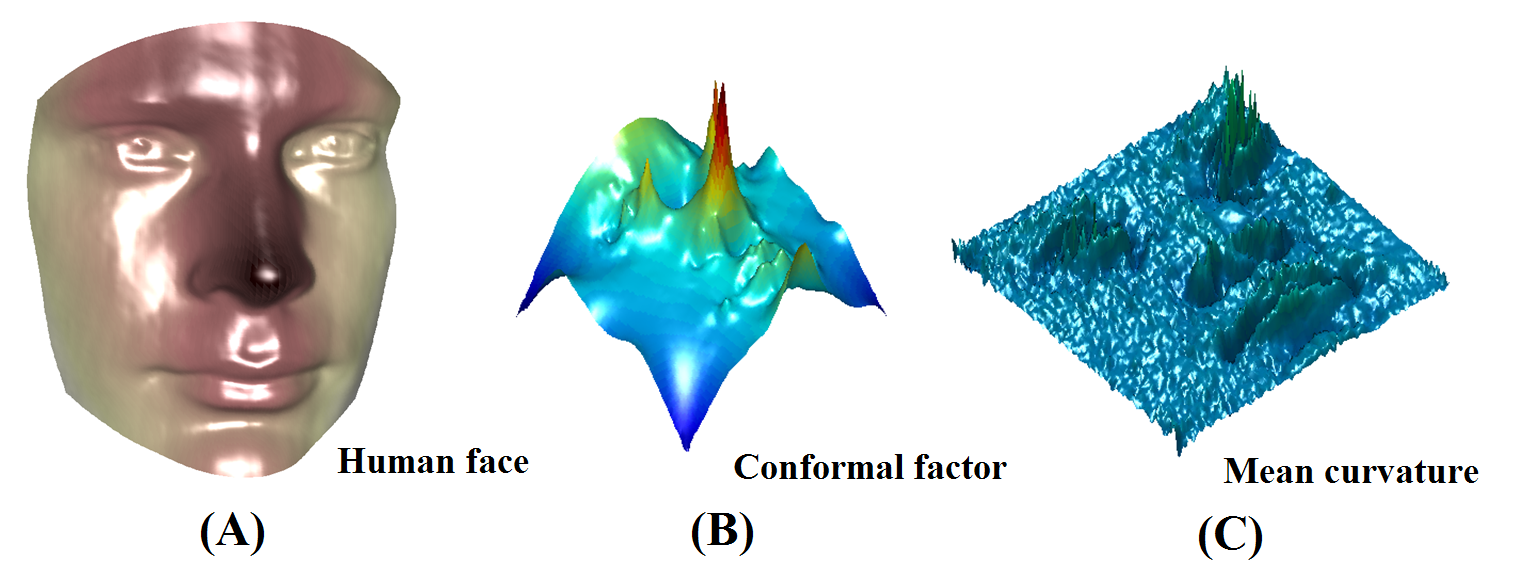}
\caption{$\lambda$-$H$ representation of the human face. (A) shows a human face. (B) shows the conformal factor $\lambda$ of the human face. (C) shows the mean curvature of the human face. The scalar functions conformal factor and mean curvature determine the surface up to a rigid motion. \label{fig:facelambdaH}}
\end{figure}

\begin{figure}[t]
\centering
\includegraphics[height=1.5in]{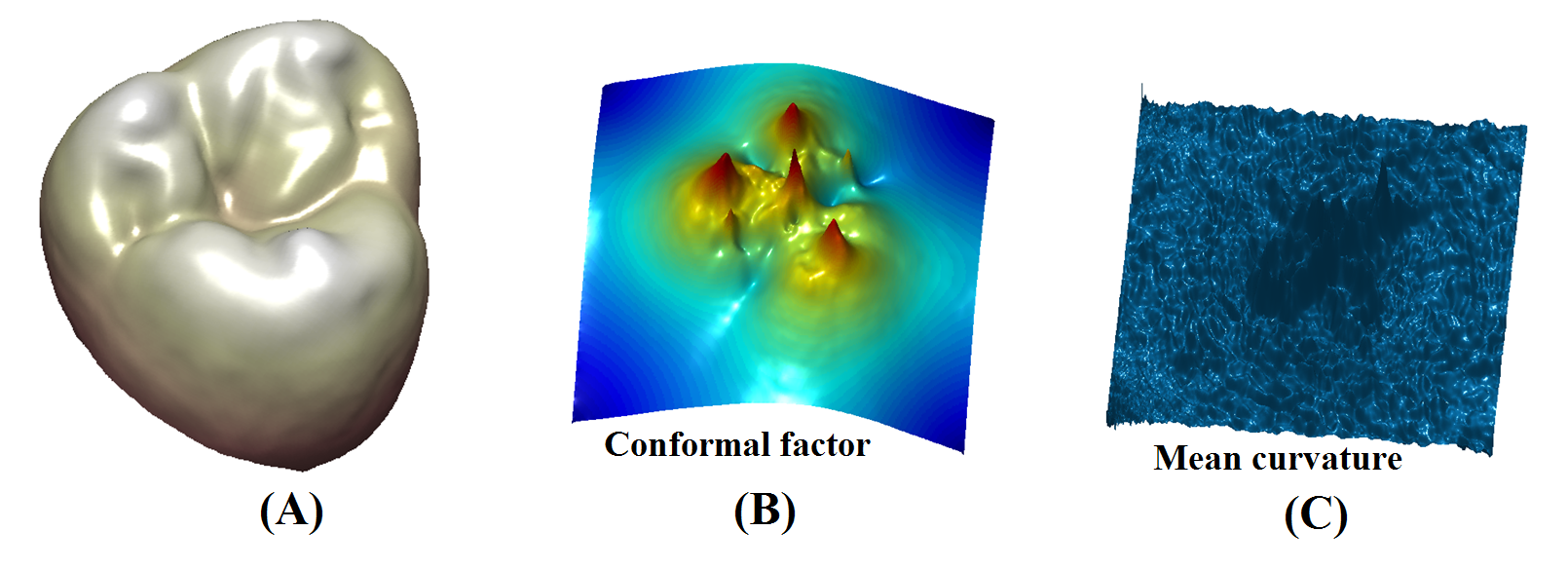}
\caption{$\lambda$-$H$ representation of a human tooth. (A) shows a human tooth. (B) shows the conformal factor $\lambda$ of the tooth. (C) shows the mean curvature of the tooth. The scalar functions conformal factor and mean curvature determine the surface up to a rigid motion. \label{fig:teethlambdaH}}
\end{figure}

\begin{figure}[t]
\centering
\includegraphics[height=1.35in]{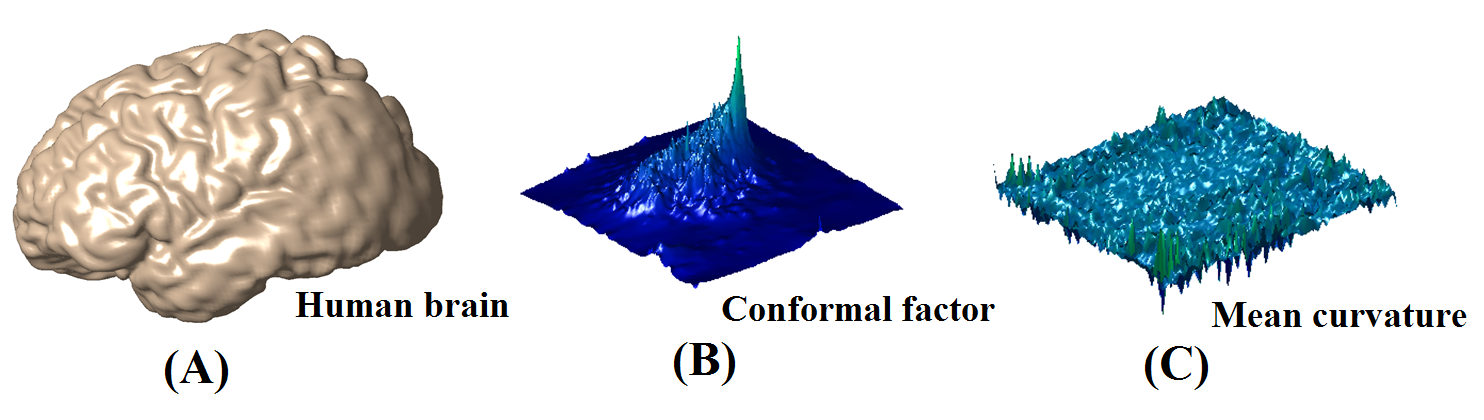}
\caption{$\lambda$-$H$ representation of the human brain surface. (A) shows a human brain surface. (B) shows the conformal factor $\lambda$ of the human brain surface. (C) shows the mean curvature of the human brain surface. \label{fig:brainlambdaH}}
\end{figure}
The above discussion can be summarized by the following theorem.

\begin{theorem}[$\lambda$-$H$ representation] Let $S$ be an open Riemann surface with mean curvature $H$. Let $\phi: D\subset \mathbb{R}^2$ be the conformal parameterization of $S$ with conformal factor $\lambda$. Suppose $U$ is a (small) boundary neighborhood of $S$. Given the boundary value $\phi|_U$ of $\phi$, $S$ is uniquely determined by $\lambda$ and $H$.
\end{theorem}

\section{Surface inpainting} In this section, we will describe a surface inpainting algorithm based on the $\lambda$-$H$ representation of the surface $S$.

Suppose $S$ is an incomplete 3D surface model with holes. Let $\Omega$ be the occluded region. Our goal is to obtain a restored surface $S'$ with holes filled based on the surface geometry, so that the basic geometry of the surface can be smoothly patched. Following conventional image inpainting algorithms, we first take an initial surface $S_0$ with the holes roughly filled. Computationally, it can be done by creating a patch piecewise linearly with boundary vertices given by $\partial D$. $S_0$ is then iteratively modified to obtain a satisfactory inpainted surface $S'$.

One way to modify $S_0$ is to update the position vector of $S_0$ iteratively. However, this method does not consider the essential geometric information of the surface. We therefore propose to inpaint the surface by considering the conformal factor $\lambda$ and mean curvature $H$, which are the two essential geometric quantities determining a Riemann surface.

Let $\phi_0: D_0 \subset \mathbb{R}^2 \to S_0$ be the conformal parameterization of $S_0$. We can compute the conformal factor $\lambda_0$ of $\phi_0$ which is given by:
\begin{equation}
\lambda_0 = <\phi_u,\phi_u>
\end{equation}

Suppose $H_0$ be the mean curvature of $S_0$. Then, ($\lambda_0$, $H_0$) determines $S_0$. As expected, the scalar functions $\lambda_0$ and $H_0$ are not smooth on the occluded region $\phi_0^{-1}(\Omega)$ in the conformal parameter domain. Our goal is to smooth out (or inpaint) $\lambda_0$ and $H_0$ on the occluded region $\phi_0^{-1}(\Omega)$, while preserving the value of them on the non-occluded region $D_0 \setminus \phi_0^{-1}(\Omega)$.

The conformal factor and mean curvature are smooth scalar functions defined on the 2D conformal parameter domain. To ensure the smoothness of the inpainted mean curvature $H'$ and conformal factor $\lambda'$, the following inpainting model is applied:

\begin{equation}\label{lambdainpaint}
\lambda' = \mathbf{argmin}_{\lambda} E(\lambda) := \mathbf{argmin}_{\lambda} \int_{\phi_0^{-1}(\Omega)} |\Delta \lambda|^2
\end{equation}

\noindent subject to the constraint that:

\begin{equation}
\lambda'|_{D_0\setminus \phi_0^{-1}(\Omega)} = \lambda_0|_{D_0\setminus \phi_0^{-1}(\Omega)}
\end{equation}

Also, we apply the same inpainting model for $H'$:

\begin{equation}\label{Hinpaint}
H' = \mathbf{argmin}_{H} E(H) := \mathbf{argmin}_{H} \int_{\phi_0^{-1}(\Omega)} |\Delta H|^2
\end{equation}

\noindent subject to the constraint that:

\begin{equation}
H'|_{D_0\setminus \phi_0^{-1}(\Omega)} = H_0|_{D_0\setminus \phi_0^{-1}(\Omega)}
\end{equation}

By minimizing the variational models (\ref{lambdainpaint}) and (\ref{Hinpaint}), we obtain an optimal (inpainted) conformal factor $\lambda'$ and $H'$. The surface associated with the pair $(\lambda', H')$ can then be reconstructed, as described in the last subsection, to obtain an inpainted surface $S'$. The surface holes will be effectively filled, with the basic geometry of the surface smoothly patched.

The minimization of the variational problems  (\ref{lambdainpaint}) and (\ref{Hinpaint}) can be done iteratively as follows:

\begin{equation}
\frac{d \lambda(t)}{dt} = - (\Delta)^2 \lambda(t) \mathrm{\ and\ } \frac{d H(t)}{dt} = - (\Delta)^2 H(t)
\end{equation}

\noindent subject to $\lambda(t)|_{D_0\setminus \phi_0^{-1}(\Omega)} = {\lambda}_0|_{D_0\setminus \phi_0^{-1}(\Omega)}$ and $H(t)|_{D_0\setminus \phi_0^{-1}(\Omega)} = H_0|_{D_0\setminus \phi_0^{-1}(\Omega)}$.

\section{Numerical algorithms}
In this section, we will describe how the proposed algorithms can be implemented. We will firstly describe how the $\lambda$-$H$ representation can be implemented. We will then explain how the surface inpainting algorithm can be implemented.

\subsection{Numerical implementation of $\lambda$-$H$ representation} Given a surface mesh $S$, the $\lambda$-$H$ representation of $S$ can be computed easily using the conformal parameterization of $S$. Many conformal parameterization algorithms have recently been developed and some of them are available online. We will not explain the conformal parameterizations in this paper, but refer the readers to \cite{Gu1,Gu2,Gu3,LuiBHF,LuiBHFHP,LuiBeltramirepresentation} for details. In this work, we use the algorithm proposed in \cite{LuiBeltramirepresentation} to obtain the conformal parameterization of $S$ onto the rectangular conformal parameter domain $R$.

\begin{figure}[t]
\centering
\includegraphics[height=1.65in]{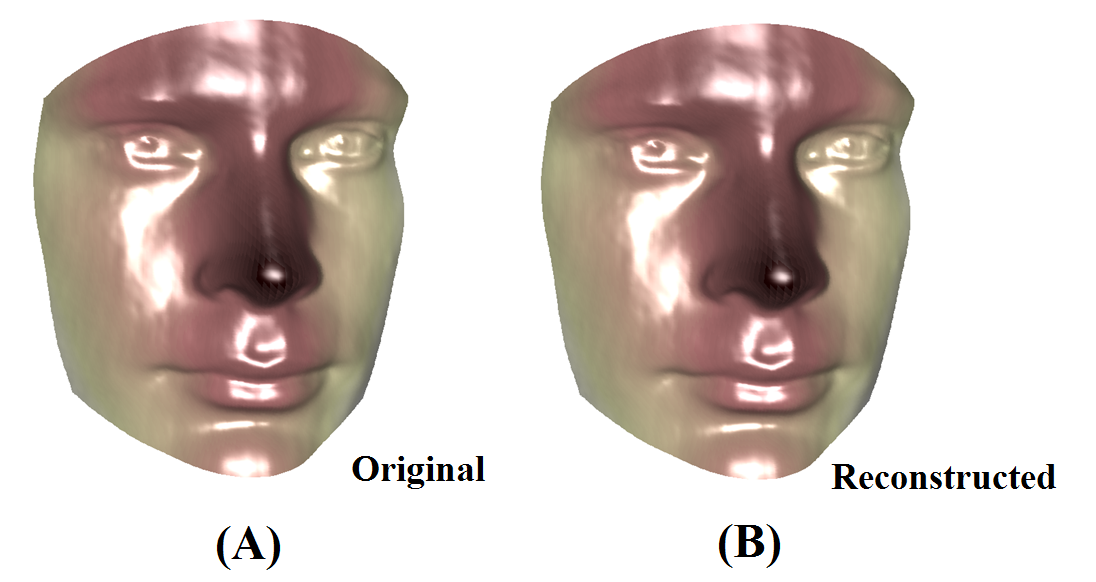}
\caption{Surface reconstruction of the human face from its $\lambda$-$H$ representation. (A) shows the original human face. The reconstructed human face from its $\lambda$-$H$ representation is shown in (B). \label{fig:facereconstruction}}
\end{figure}

\begin{figure}[t]
\centering
\includegraphics[height=1.65in]{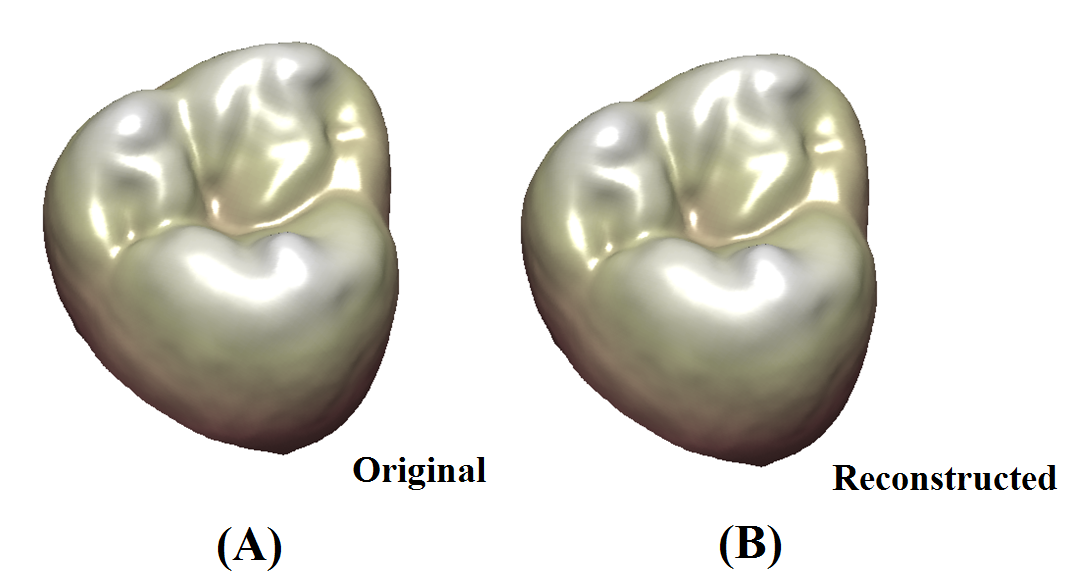}
\caption{Surface reconstruction of the human tooth from its $\lambda$-$H$ representation.  (A) shows the original tooth. The reconstructed tooth surface from its $\lambda$-$H$ representation is shown in (B). \label{fig:teethreconstruction}}
\end{figure}

\begin{figure}[t]
\centering
\includegraphics[height=1.5in]{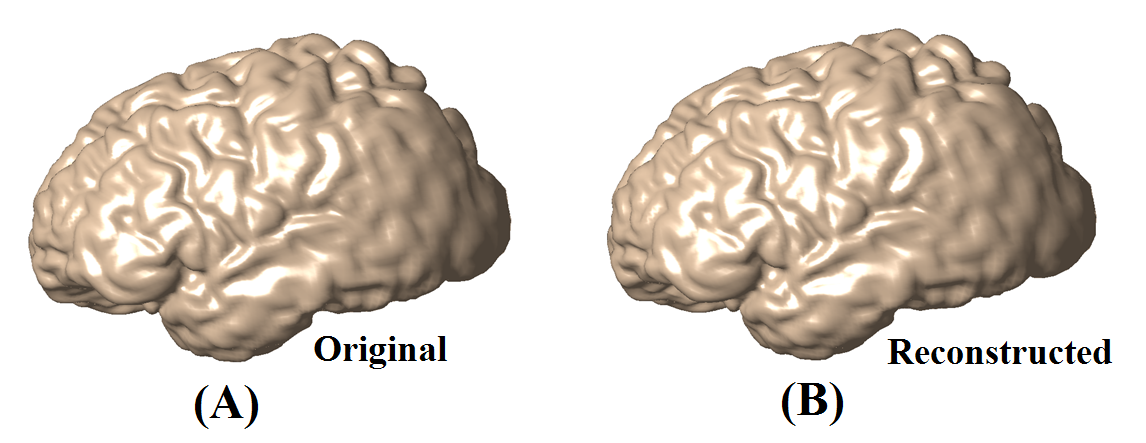}
\caption{Surface reconstruction of the human brain from its $\lambda$-$H$ representation.  (A) shows the original brain surface. The reconstructed brain surface from its $\lambda$-$H$ representation is shown in (B). \label{fig:brainreconstruction}}
\end{figure}

Let $\vec{\phi}(u,v) := (X(u,v),Y(u,v),Z(u,v))$ be the conformal parameterization of the surface mesh $S$. We consider the finite difference discretization of the partial derivatives on the rectangular conformal parameter domain $R : = [0,1]\times [0,K]$. Let $u_i = i/N$, $v_j = Kj/N$ and $h=1/N$. We denote $\vec{\phi}_{ij} = \vec{\phi} (u_i,v_j);$  $H_{ij} = H(u_i,v_j);$ and $\lambda_{ij} = \lambda(u_i,v_j)$. Using equation (\ref{I}), the conformal factor $\lambda_{ij}$ can be computed by:

\begin{equation}
\lambda_{ij}^2 = [\ |\frac{\vec{\phi}_{i+1,j} - \vec{\phi}_{i-1,j}}{2h}|^2 + |\frac{\vec{\phi}_{i,j+1} - \vec{\phi}_{i,j-1}}{2h}|^2\ ]/2
\end{equation}

And using equation (\ref{H}), the mean curvature $H_{ij}$ can be computed by:

\begin{equation}
H_{ij} = \frac{1}{2\lambda_{ij}}\mathbf{sign}(\vec{\phi}_{ij})|\Delta \vec{\phi}_{ij}|
\end{equation}

\noindent where:

\begin{equation}
\begin{split}
& \Delta \vec{\phi}_{ij} = \frac{\vec{\phi}_{i+1,j}+\vec{\phi}_{i-1,j}- 4\vec{\phi}_{i,j} + \vec{\phi}_{i,j+1}+\vec{\phi}_{i,j-1}}{h^2}\\
&\mathbf{sign}(\vec{\phi}_{ij}) = \mathbf{sign}(<\Delta \vec{\phi}_{ij}, \vec{n}_{ij}>)\\
& \vec{n}_{ij} = (\frac{\vec{\phi}_{i+1,j} - \vec{\phi}_{i-1,j}}{2h})\times (\frac{\vec{\phi}_{i,j+1} - \vec{\phi}_{i,j-1}}{2h})/||(\frac{\vec{\phi}_{i+1,j} - \vec{\phi}_{i-1,j}}{2h})\times (\frac{\vec{\phi}_{i,j+1} - \vec{\phi}_{i,j-1}}{2h})||
\end{split}
\end{equation}

Conversely, given the $\lambda$-$H$ representation of $S$, we can reconstruct the surface by solving the equation (\ref{frame}). Let $(U_{ij}, V_{ij},W_{ij}) = ((\vec{\phi}_z)_{ij}, (\vec{\phi}_{\bar{z}})_{ij}, \vec{n}_{ij})$ be the natural frame at $(u_i,v_j)$. Consider the central difference discretization of $\frac{\partial}{\partial z}$ and $\frac{\partial}{\partial \bar{z}}$ as follows:

\begin{equation}
\begin{split}
\frac{\partial f}{\partial z} & := \frac{1}{2}[\left( \frac{f_{i+1,j}- f_{i-1,j}}{2h}\right) - \sqrt{-1} \left( \frac{f_{i,j+1}- f_{i,j-1}}{2h}\right)]\\
\frac{\partial f}{\partial \bar{z}} & := \frac{1}{2}[\left( \frac{f_{i+1,j}- f_{i-1,j}}{2h}\right) + \sqrt{-1} \left( \frac{f_{i,j+1}- f_{i,j-1}}{2h}\right)]
\end{split}
\end{equation}

Hence, $(\lambda_z){ij}$ can be computed for all $0 \leq i,j \leq N$. Together with the boundary conditions, the natural frame $(U_{ij}, V_{ij},W_{ij})$ can then be computed by solving the following linear system, which approximates the solution of the equation (\ref{frame}):

\begin{equation}
\frac{1}{2h}\left(\begin{array}{c}
(U_{i+1,j}- U_{i-1,j}) - \sqrt{-1}(U_{i,j+1}- U_{i,j-1})\\
(V_{i+1,j}- V_{i-1,j}) - \sqrt{-1}(V_{i,j+1}- V_{i,j-1})\\
2h W_{ij}
\end{array}\right) = \left(\begin{array}{c}
\frac{2}{\lambda_{ij}} (\lambda_z)_{ij}U_{ij} + \mu_{ij}W_{ij}\\
\frac{\lambda_{ij}^2}{2} H_{ij} W_{ij}\\
-H_{ij} U_{ij} - \frac{2\mu_{ij} V_{ij}}{\lambda_{ij}^2}
\end{array}\right)
\end{equation}

Once the natural frame is obtained, the position vector $\vec{\phi}_{ij}$ (and hence the surface mesh $S$) can be reconstructed by solving the linear system (with the given boundary conditions):

\begin{equation}
\frac{1}{2h}\left(\begin{array}{c}
\vec{\phi}_{i+1,j} - \vec{\phi}_{i-1,j}\\
\vec{\phi}_{i,j+1} - \vec{\phi}_{i,j-1}
\end{array}\right) = \left(\begin{array}{c}
(U_{ij} + V_{ij})/2\\
\sqrt{-1}(U_{ij} - V_{ij})/2
\end{array}\right)
\end{equation}

\subsection{Numerical implementation of surface inpainting}
Given a surface mesh with holes $S$, we first obtain an initial mesh $S_0$ with the holes roughly filled. It can be obtained by the delaunay triangulation of the boundary vertices of the surface holes. The $\lambda$-$H$ representation of $S_0$ can then be computed, which are two non-smooth scalar functions of the conformal factor $\lambda^0$ and mean curvature $H^0$. We then inpaint (or smooth out) $\lambda_0$ and $H_0$ at the occluded regions on the conformal parameter domain, using the inpainting models (\ref{lambdainpaint}) and (\ref{Hinpaint}).

Consider the discretization of the Laplacian on the conformal parameter domain as follows:

\begin{equation}\nonumber
(\Delta f)_{ij} = \frac{f_{i+1,j} + f_{i-1,j} - 4f_{ij} + f_{i,j+1} + f_{i,j-1}}{h^2}
\end{equation}

With this discretization, the inpainting models on $\lambda^0$ and $H^0$ can be implemented as follows:

\begin{equation}
\begin{split}
\lambda^n_{ij} &= \lambda^{n-1}_{ij} + dt [\Delta (\Delta \lambda^{n-1})]_{ij}\\
H^n_{ij} &= H^{n-1}_{ij} + dt [\Delta (\Delta H^{n-1})]_{ij}
\end{split}
\end{equation}

\noindent with the constraints that $\lambda^n_{ij} = \lambda^0_{ij}$ and $H^n_{ij} = H^0_{ij}$ for all $(u_i,v_j)$ that are not in the occluded regions. $dt$ is the time step.

\begin{figure}[t]
\centering
\includegraphics[height=3.8in]{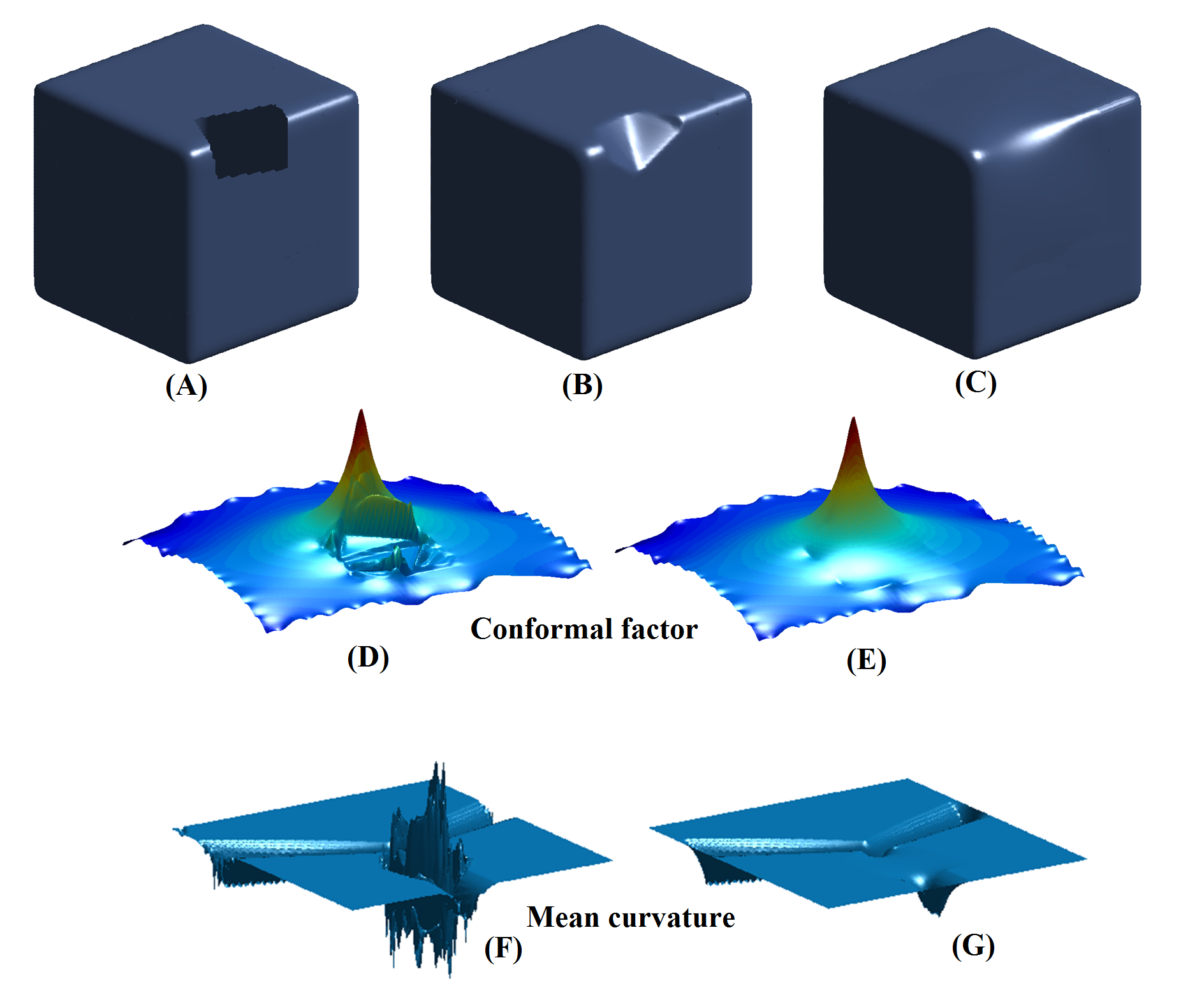}
\caption{Surface inpainting of the cubic surface. (A) shows the cubic surface with an occlusion. (B) shows the inpainted surface obtained from the conventional smooth hole-filling method. (C) shows the inpainted surface using our proposed conformal approach.  (D) and (E) shows the initial conformal factor and inpainted (smoothed) conformal factor respectively. (F) and (G) shows the initial mean curvature and inpainted (smoothed) mean curvature respectively. \label{fig:cube}}
\end{figure}

\begin{figure}[t]
\centering
\includegraphics[height=4.45in]{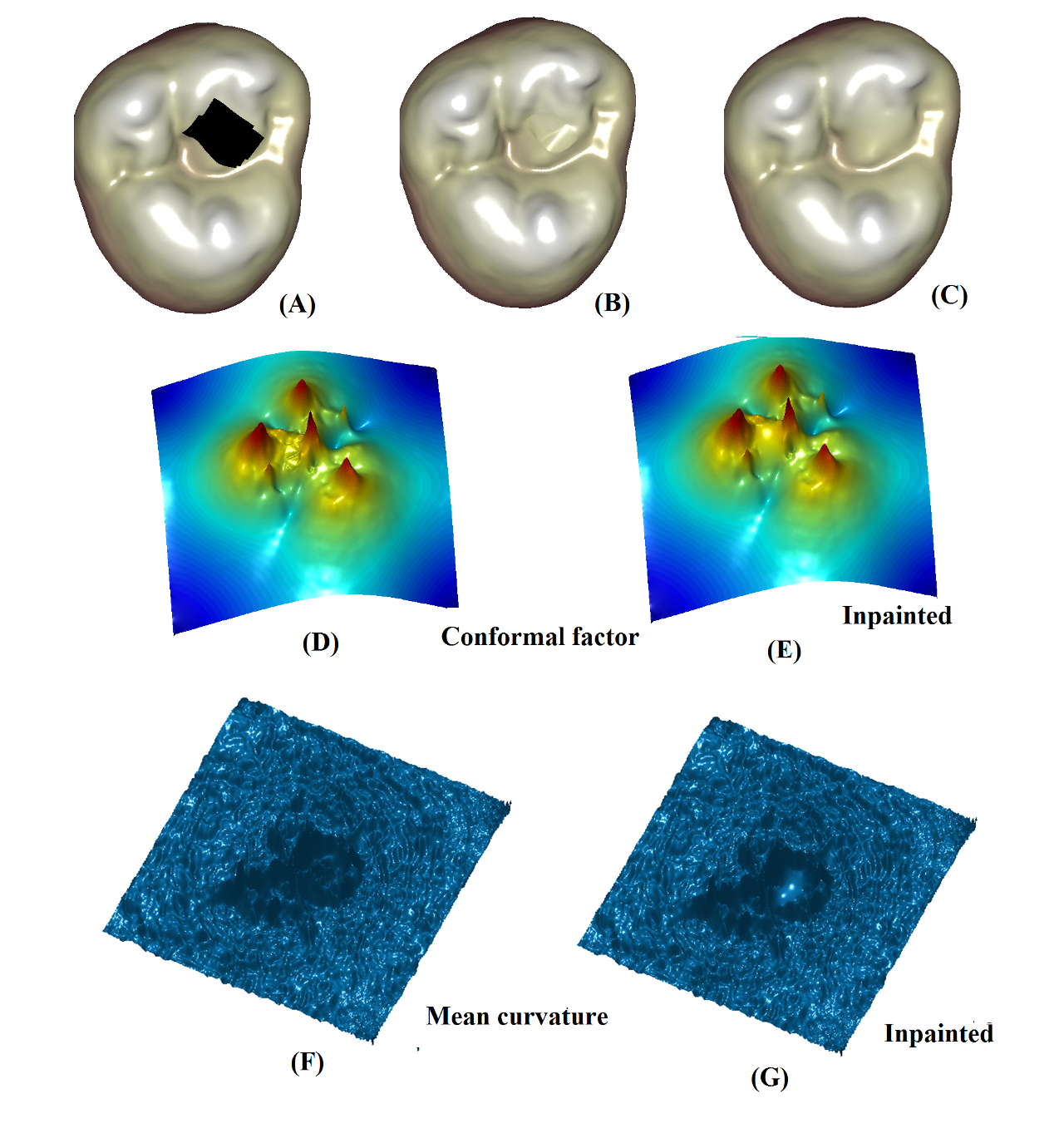}
\caption{Surface inpainting of the human tooth surface. (A) shows the tooth surface with an occlusion. (B) shows the inpainted surface obtained from the conventional smooth hole-filling method. (C) shows the inpainted surface using our proposed conformal approach.  (D) and (E) shows the initial conformal factor and inpainted (smoothed) conformal factor respectively. (F) and (G) shows the initial mean curvature and inpainted (smoothed) mean curvature respectively. \label{fig:teeth}}
\end{figure}

\begin{figure}[t]
\centering
\includegraphics[height=5.15in]{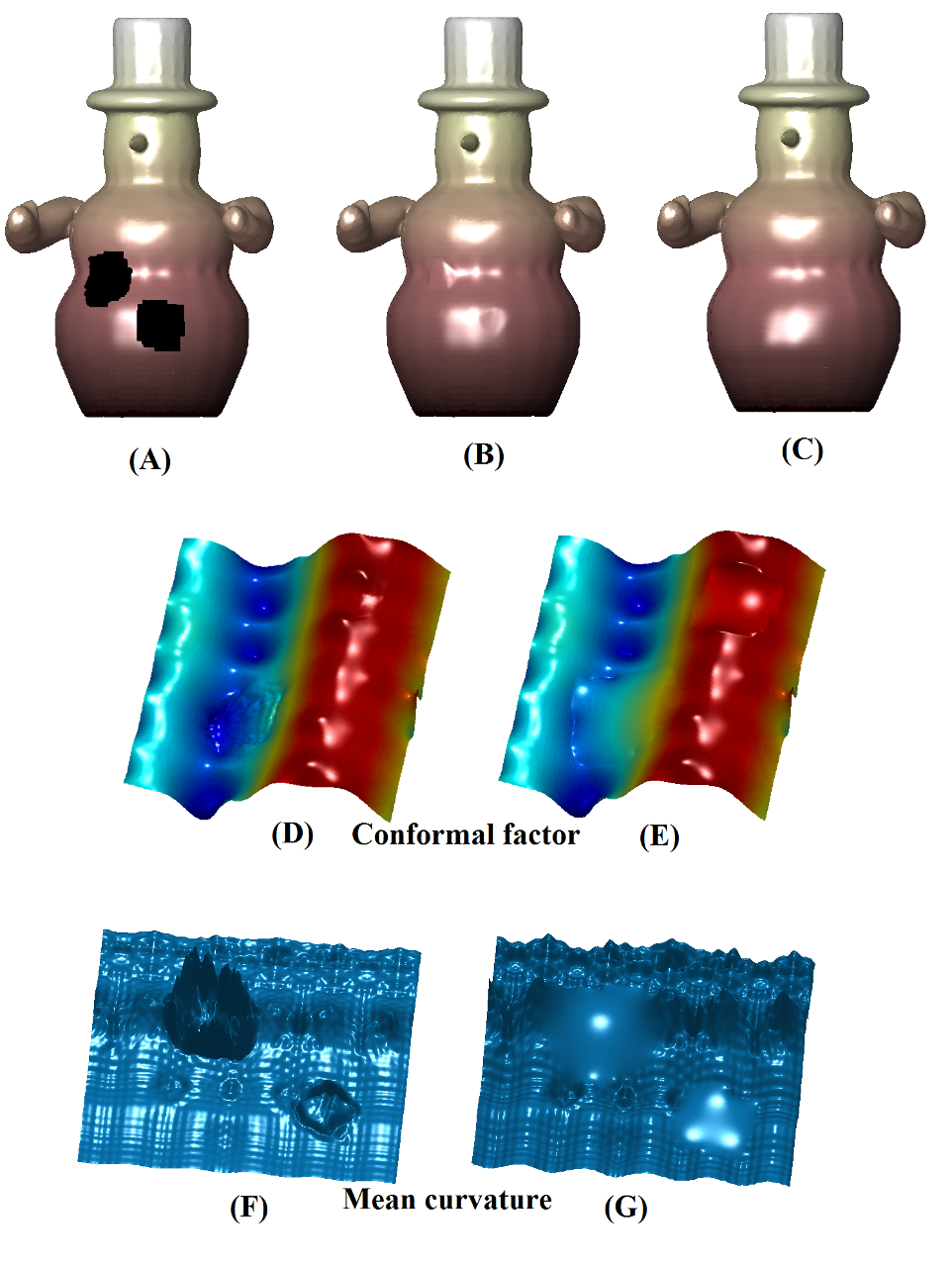}
\caption{Surface inpainting of the snowman surface. (A) shows the snowman surface with two holes. (B) shows the inpainted surface obtained from the conventional smooth hole-filling method. (C) shows the inpainted surface using our proposed conformal approach.  (D) and (E) shows the initial conformal factor and inpainted (smoothed) conformal factor respectively. (F) and (G) shows the initial mean curvature and inpainted (smoothed) mean curvature respectively.  \label{fig:snowman}}
\end{figure}

\begin{figure}[t]
\centering
\includegraphics[height=4.15in]{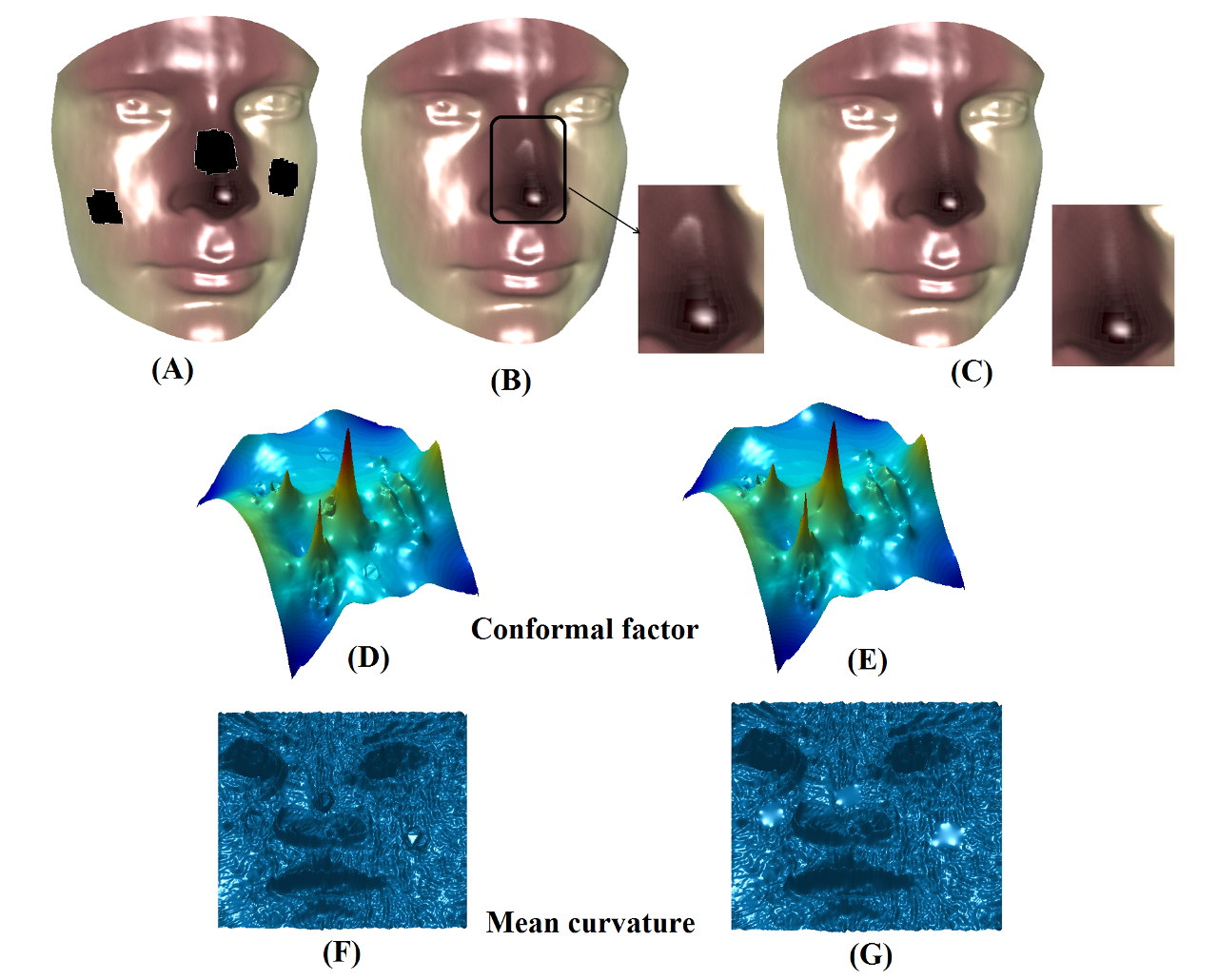}
\caption{Surface inpainting of the human face surface. (A) shows the human face with three holes. (B) shows the inpainted surface obtained from the conventional smooth hole-filling method. (C) shows the inpainted surface using our proposed conformal approach.  (D) and (E) shows the initial conformal factor and inpainted (smoothed) conformal factor respectively. (F) and (G) shows the initial mean curvature and inpainted (smoothed) mean curvature respectively. \label{fig:faceinpaint}}
\end{figure}

\begin{figure}[t]
\centering
\includegraphics[height=4.75in]{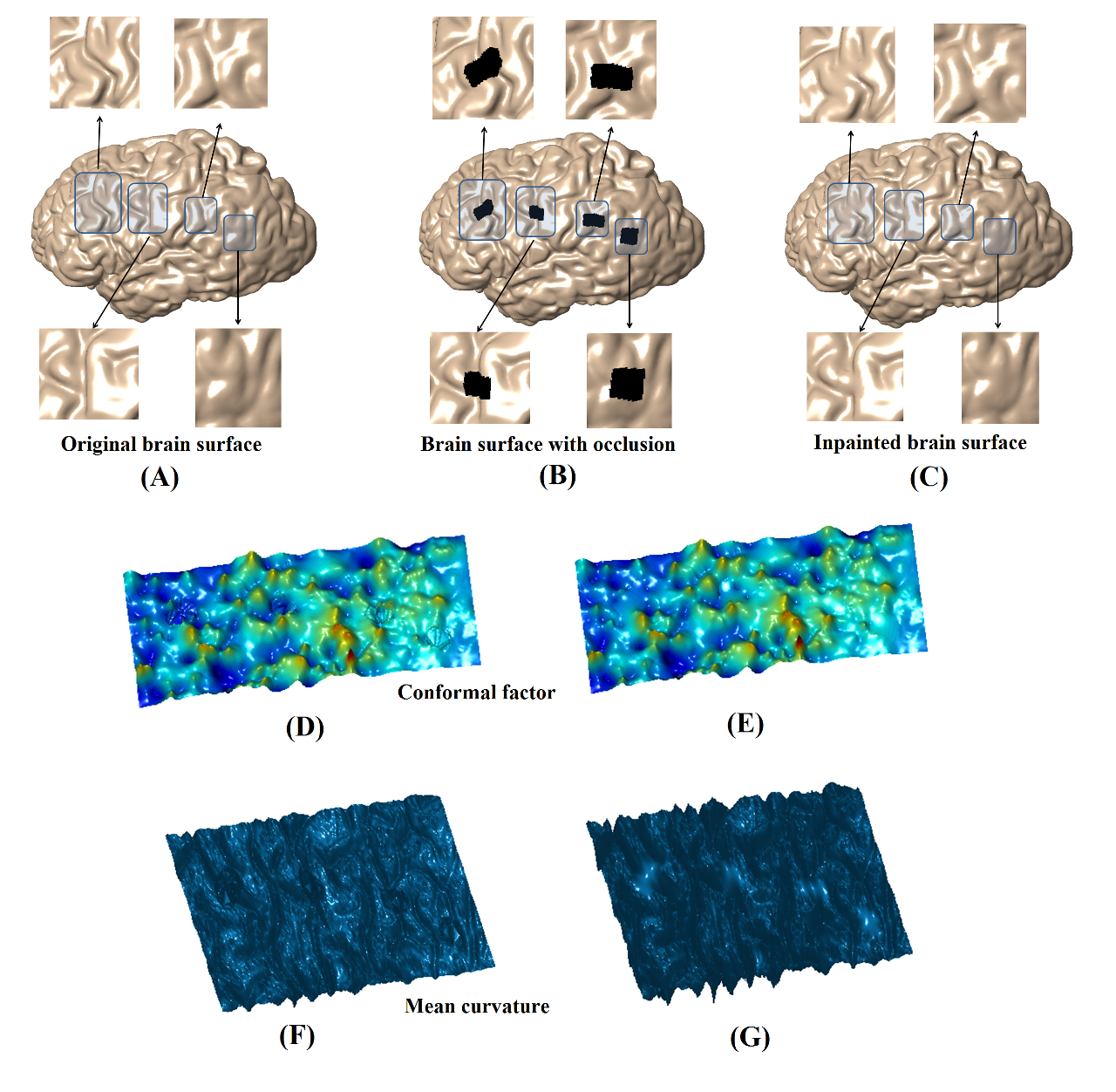}
\caption{Surface inpainting of the human brain surface. (A) shows the original brain surface. We remove four regions on the brain surface, as shown in (B). (C) shows the inpainted surface using our proposed conformal approach. Note that the inpainted patches closely resemble to the original one. (D) and (E) shows the initial conformal factor and inpainted (smoothed) conformal factor respectively. (F) and (G) shows the initial mean curvature and inpainted (smoothed) mean curvature respectively.  \label{fig:braininpaint}}
\end{figure}

\section{Experimental results}
To examine the effectiveness of the proposed method, we have tested our algorithms on synthetic data together with real surface data. In this section, we will describe the experimental results in details.

\subsection{Surface reconstruction from $\lambda$-$H$ representation}
We first examine how well the $\lambda$-$H$ representation can represent a Riemann surface. In Figure \ref{fig:facelambdaH}, we compute the conformal factor $\lambda$ and mean curvature $H$ of a human face. This gives the $\lambda$-$H$ representation of the human face. Using $\lambda$-$H$ representation, we can accurately reconstruct the original surface. The reconstructed surface is shown in Figure \ref{fig:facereconstruction}(B), which closely resembles to the original surface as shown in Figure \ref{fig:facereconstruction}(A).

We also compute the $\lambda$-$H$ representation of a human tooth, as shown in Figure \ref{fig:teethlambdaH}. From the $\lambda$-$H$ representation, the tooth surface can be effectively and accurately reconstructed. Figure \ref{fig:teethreconstruction}(B), shows the reconstructed tooth surface from its $\lambda$-$H$ representation, which again closely resembles to the original tooth surface as shown in Figure \ref{fig:teethreconstruction}(A).

Finally, we compute the $\lambda$-$H$ representation of a human brain, whose geometry is complicated. The conformal factor and mean curvature of the brain surface is shown in Figure \ref{fig:brainlambdaH}. From the $\lambda$-$H$ representation, the brain surface can be reconstructed as shown in Figure \ref{fig:brainreconstruction}. The reconstructed surface is close to the original one. It demonstrates that the proposed geometric representation is effective, even for surfaces with complicated geometry.

\subsection{Surface inpainting}
Next, we test the surface inpainting algorithm on synthetic data and real surface data. In Figure \ref{fig:cube}, we test the algorithm on a cubic surface with a hole. To inpaint the surface, a conventional way is to smoothly fill in the surface hole with its boundary representation in 3D \cite{InpaintingClarenz,InpaintingDavis}. Although the surface can be smoothly filled, the inpainted surface cannot follow the surface geometry. As shown in (B), the inpainted surface generated by the conventional smooth hole-filling algorithm cannot preserve the sharp edge of the cubic surface. The sharp edge is cropped off. (C) shows the inpainted surface using our proposed conformal approach. The inpainted surface can better preserve the sharp edge. (D) and (E) shows the initial conformal factor and inpainted (smoothed) conformal factor respectively. (F) and (G) shows the initial mean curvature and inpainted (smoothed) mean curvature respectively. By inpainting the mean curvature and conformal factor, the inpainted surface follows the surface geometric pattern of the non-occluded regions of the surface.

In Figure \ref{fig:teeth}, we examine the proposed inpainting algorithm on a human tooth surface. (A) shows a tooth surface with an occlusion. (B) shows the inpainted surface obtained from the conventional smooth hole-filling algorithm. Again, the inpainted surface does not follow the surface geometry of the non-occluded regions of the surface. (C) shows the inpainted surface generated by our proposed inpainting method. The inpainted surface is more natural which follows the surface geometry. The basic geometry can be smoothly patched. (D) and (E) shows the initial conformal factor and inpainted (smoothed) conformal factor respectively. (F) and (G) shows the initial mean curvature and inpainted (smoothed) mean curvature respectively.

In Figure \ref{fig:snowman}, we show the surface inpainting result of a snowman surface. (A) shows the snowman surface with 2 holes. (B) shows the inpainted surface obtained from the conventional smooth hole-filling algorithm. The inpainted surface patches are flatten and the inpaint surface does not follow the surface geometry. (C) shows the inpainted surface using our method. The inpainted surface is more natural which follows the surface geometry.

Figure \ref{fig:faceinpaint} shows the example of surface inpainting of a human face. (A) shows a human face with 3 holes, one on the nose and two on the cheeks. (B) shows the inpainted surface obtained from the conventional smooth hole-filling algorithm. The inpainted surface is not natural and does not following the surface geometry. For example, the inpainted surface patch (inside the box) for the hole on the nose is flatten. (C) shows the inpainting result obtained by our algorithm. The reconstruct surface is more natural and follows the surface geometric pattern.

Finally, in Figure \ref{fig:braininpaint}, we test the inpainting algorithm on a brain surface, which has more complicated geometry. (A) shows the original brain surface. We remove four regions on the brain surface, as shown in (B). The inpainted surface generated by our algorithm is shown in (C). Note that the inpainted surface patches closely resemble to the original ones (compare with (A)). This example demonstrates that our proposed inpainting algorithm perform well even on surfaces with complicated geometry.

\section{Conclusion and future works}
In this work, we propose a novel conformal approach for surface inpainting to fill in missing holes on an incomplete 3D surface model. It has important applications in different fields, such as in medical imaging, computer graphics and computer visions. In order to fully manipulate the surface geometry of the captured 3D model, we consider in this work a new representation of 3D surfaces based on their conformal factor $\lambda$ and mean curvature $H$. According to Riemann surface theory, a Riemann surface can be uniquely determined by its conformal factor and mean curvature up to a rigid motion. Given a Riemann surface $S$, its mean curvature $H$ and conformal factor $\lambda$ can be computed easily through its conformal parameterization. Conversely, given $\lambda$ and $H$, a Riemann surface can be uniquely reconstructed by solving the Gauss-Codazzi equation on the conformal parameter domain. Hence, the conformal factor and the mean curvature are two geometric quantities fully describing the surface. In this paper, we develop a numerical algorithm to reconstruct surfaces from their conformal factor and mean curvature by solving systems of linear equations.  Given the $\lambda$-$H$ representation of an incomplete surface, the problem of surface inpainting can be reformulated into the image inpainting problem of $\lambda$ and $H$ on the conformal parameter domain. Using the conventional image inpainting techniques, both $\lambda$ and $H$ can be inpainted on the parameter domain. A Riemann surface can then be reconstructed which restores the 3D surface with surface holes.  Since the inpainting model is based on the geometric quantities $\lambda$ and $H$, the restored surface follows the surface geometric pattern as much as possible. We test the proposed algorithm on synthetic data, 3D human face data and MRI-derived brain surfaces. Experimental results show that our algorithm can effectively inpaint surface holes based on the geometry of the surface to restore the incomplete 3D surface models.

Besides, we observe from our experiments that the $\lambda$-$H$ representation of a surface is a useful geometric representation of a surface mesh capturing the most essential geometric information. By manipulating $\lambda$ and $H$, the surface mesh can be edited based on the surface geometry. Thus, we believe there are several potential applications of this representation to other surface mesh editing problems, such as surface mesh compression, surface mesh denoising and so on. In the future, we will consider applying the $\lambda$-$H$ representation for surface stitching. We will also consider extending the $\lambda$-$H$ representation to point cloud. With that, point cloud compression, point cloud inpainting and point cloud denoising can be easily achieved.
%%%%%%%%%%%%%%%%%%%%%%%%%%%%%%%%%%%%%%%%%%%%%%%%%%%%%%%%%%%%%%%%%%%%%
% REFERENCES
%%%%%%%%%%%%%%%%%%%%%%%%%%%%%%%%%%%%%%%%%%%%%%%%%%%%%%%%%%%%%%%%%%%%%
%{
%\bibliographystyle{ieee}
%\bibliography{egbib}
%}

\end{document}